# Fine-tuning the Ant Colony System algorithm through Particle Swarm Optimization.

D. Gómez-Cabrero[+], D. N. Ranasinghe*

*University of Colombo School of Computing
35, Reid Avenue, Colombo 7, Sri Lanka
Email: dnr@ucsc.cmb.ac.lk
[+]Departamento de Estadística e Investigación Operativa,
Universitat de Valencia, Calle Doctor Moliner s/n, Burjassot, España
Email: lunacab@yahoo.es

*Abstract*—Ant Colony System (ACS) is a distributed (agent-based) algorithm which has been widely studied on the Symmetric Travelling Salesman Problem (TSP). The optimum parameters for this algorithm have to be found by trial and error. We use a Particle Swarm Optimization algorithm (PSO) to optimize the ACS parameters working in a designed subset of TSP instances. First goal is to perform the hybrid PSO-ACS algorithm on a single instance to find the optimum parameters and optimum solutions for the instance. Second goal is to analyze those sets of optimum parameters, in relation to instance characteristics. Computational results have shown good quality solutions for single instances though with high computational times, and that there may be sets of parameters that work optimally for a majority of instances.

## I. Introduction

Heuristics algorithms, and their higher forms meta-heuristics, have been widely used to find reasonable good solutions for NP-Hard problems. The performance of those heuristics is based on the optimum set-up of a set of parameters. The problem of fine-tuning those parameters is also a hard problem.

Fine-Tuning is an unavoidable task that needs a correct design of experiments. For statistical experimental design and analysis for heuristics see [4]. Heuristics as CALIBRA [1] and F-Race [3] have been designed as procedures for finding optimal parameters. Fine-tuning can be understood as "to find the set of parameters performing well on a wide set of instances".

We have designed an algorithm, using a Particle Swarm Optimization (PSO) framework, to optimize the parameters of the ACS algorithm working on a "single" Symmetric Travelling Salesman Problem (TSP) instance. For each instance the algorithm computes an "optimal" set of ACS parameters. A second goal is to analyze jointly those sets of parameters, their performance on all instances (not only their related instance) and the characteristics of their related instances for the purpose of finding correlations.

Ant Colony Optimization (ACO) is a generic framework for ant-based optimization heuristic algorithms. In ACO algorithms ants are agents which, in the TSP case, construct tours by moving from city to city on the graph problem. Those ants are sharing information using a pheromone trail.

The first ACO algorithm, called Ant-System, was proposed in [5], [7] and [8]. A full review of ACO algorithms and applications can be found in [9]. ACS is a version of the Ant System that modifies the updating of the pheromone trail, see [6] and [11]. We have chosen this ACS algorithm to work with because of the theoretical background we have found on it, see [6] and [9], and the previous fine-tuning research on the parameters by [16].

PSO is a Swarm Intelligence method for global optimization. Given a domain $D$, there is defined a function $f : D \rightarrow \Re$ assigning to each point of $D$ a fitness value. In PSO there is a population (swarm) of individuals, named particles, moving on the domain and adjusting their trajectory to its own previous best position and the previous best position of the neighbourhood. We will use the Global PSO version which considers the neighbourhood as all the swarm. For a further review of the PSO [10], [13] and [14]. An example of PSO applications can be found in [12] and [15]. We have chosen PSO because it has an easy implementation for integer and real parameters, and, as genetic algorithms, it performs a "blind" search on all the possible sets of parameters.

In our algorithm the domain of the PSO will be all possible sets of parameters for ACS. For a position of a particle we compute the fitness by running the ACS algorithm with the parameters given by the position on a TSP instance.

Section II describes ACS and PSO algorithms. Section III describes the PSO-ACS algorithm, the parameters used on the PSO and the initial population of particles (set of feasible parameters for the ACS algorithm). Computational results are given in section IV and finally in section V conclusions are set.

## II. PSO AND ACS ALGORITHMS

Let us introduce some notation about a TSP graph, where for a given instance denote $V$ as the set of vertices, $E$ as the set of edges (shortest paths between the vertices), $c_e$ cost of traversing edge $e \in E$.

### A. Ant Colony System (ACS)

The ACS works as follow: it has a population of $n$ ants. Let denote for each arc $e=(i,j)$ in the TSP-instance graph an initial heuristic value $\eta_e$ and an initial pheromone value $\tau_e$. $\eta_e$ is originally set to the inverse of the cost of traversing the edge $e$. $\tau_e$ is initially set to $\tau_0 = 1 \backslash L_{nn}$ for all edge $e$, where $L_{nn}$ is equal to the inverse of the tour length computed by the nearest-neighbour-heuristic algorithm.

Let $q_0, \alpha, \rho \in [0,1]$, be real values and $\varphi, \beta$ integer values between 0 and 8. For each vertex $s \in V$ a neighbour set is defined among the nearest vertices, $N(s)$. For a given ant $r$ let be $NV(r)$ the set of non-visited vertices. We denote $J_r(s)=N(s) \cap NV(r)$ the set of non-visited vertices among the neighbour set for a given vertex $s$ and a given ant $r$.

At each iteration, each ant constructs a tour solution for the TSP-instance. The constructions phase works as follows:

Each ant is initially set in a randomly vertex, then at each step the entire ants make a movement to a non-visited vertex. Given and ant $r$ with an actual position (vertex) $s$, $p_{k(rs)}$ is computed as a reference value for visiting or not vertex $k$, where:

$$p_{krs} = \begin{cases} \dfrac{[\tau_{(s,k)}]^{\varphi}[\eta_{(s,k)}]^{\beta}}{\sum_{u \in J(s)}[\tau_{(s,u)}]^{\varphi}[\eta_{(s,u)}]^{\beta}} & \text{if } k \in J(s) \\ 0 & \text{otherwise} \end{cases} \quad (1)$$

This formula includes a small modification respect to the original ACS algorithm including $\varphi$ as exponent of the pheromone level, this will allow us a deeper research on the effects of the possible combinations of $\varphi, \beta$ parameters.

A sample random value $q$ is computed. If $q \leq q_0$ we visit the city $k \in V$ with maximum $p_{k(rs)}$ (exploitation of the knowledge) otherwise ACS follows a random-proportional rule on $p_{k(rs)}$ for all $k \in V$ (biased exploration). If there are no non-visited vertex on the neighbour of vertex $s$, we extend (1) to all vertices in $NV(r)/N(s)$ (those not visited by ant $r$ and not included in the neighbour of $s$) and visits the vertex with maximum $p_{k(rs)}$.

After an arc is inserted into a route (a new vertex is visited), its pheromone trail is updated. This phase is called Local Update, and for a inserted $e \in E$:

$$\tau_e = (1-\alpha)\tau_e + \alpha\tau_0 \quad (2)$$

This reduces the pheromone level in the visited arcs and the exploration in the set of possible tours is increased.

When all the tours have been computed a Global Update phase is done. For each edge $e$ pertaining to the global-best-tour found:

$$\tau_e = (1-\alpha)\tau_e + \alpha\Delta\tau_e \quad (3)$$

$$\Delta\tau_e = 1/L_{gb} \quad (4)$$

where $L_{gb}$ is the length of the global-best-tour found.

In the original ant algorithm, and in most of the later versions, pheromone global update is performed in all the edges; ACS only updates pheromone level in the set of edges pertaining to the best tour.

We consider a trial as a performance of 1000 iterations. The lowest length tour found after all iterations are finished is the best solution found by the trial.

A feasible set of parameters for running ACS is a combination of feasible $q_0, \varphi, \beta, \alpha, \rho$, number of ants ($na$) and a concrete neighbour definition.

### B. Particle Swarm Optimization (PSO)

As described by Eberhert and Kennedy [10], [13], [14] PSO is an adaptative algorithm based on a social environment where a set of particles, called population, are visiting possible "positions" of a given domain. Each position has a fitness value (and it can be computed). At each iteration particles will move returning stochastically toward its previous best fitness position and to the population best fitness position. All the particles of the population are sharing information of the best areas to search.

Let denote $P$ as the set of parameters and let define $PO$ as the population of particles. At each iteration $x_{fp}$ and $v_{fp}$ denotes respectively the actual position and the actual velocity of parameter $p \in P$ of the particle $f \in PO$.

The movement of the particles are defined by the following equations:

$$v_{fp} = wv_{fp} + r_1 c_1 (x_{fp} - bl_{fp}) + r_2 c_2 (x_{fp} - bg_p) \quad (5)$$

$$x_{fp} = x_{fp} + \chi v_{fp} \quad (6)$$

Where $c_1, c_2$ are integer non-negative values, named cognitive and social respectively, $r_1, r_2$ are sample random values in [0,1], $w$ and $\chi$ are non-negative real values, named respectively inertia weight and constriction factor, $bg_p$ is the value of parameter p pertaining to the best set of parameters found by the population (social knowledge) and $bl_{fp}$ is the value of parameter p pertaining to the best parameters set found by particle f (self-knowledge). In (5) the first factor refers to the previous velocity, second and third factors are related respectively to the distance to the best set of parameters found by the particle and to the distance to the best set of parameters found by the population.

### III. PSO-ACO ALGORITHM

*A. PSO Parameters, initial particles and fitness value.*

The algorithm is run each time on a single TSP-instance. The set of parameters of ACS that define a point in the PSO domain are $q_0$, $\varphi$, $\beta$, $\phi$, $\alpha$, $\rho$ and the number of ants ($na$). Most of them have been previously explained in the section II. $\phi$ denotes the percentage of vertices that will be included into $N(v)$ for any vertex $v \in V$, so given a $v \in V$ and $\phi=0.5$, $|N(v)|=\lceil \phi*|V| \rceil$. The ranges of each parameter are shown in Table I, where each parameter pertains to its related ]Minimum, Maximum].

TABLE I
RANGE OF ACS PARAMETERS

| | $q_0$ | $\varphi$ | $\beta$ | $\rho$ | $\alpha$ | $\phi$ | $na$ |
|---|---|---|---|---|---|---|---|
| Minimum | 0 | -1 | -1 | 0 | 0 | 0 | 0 |
| Maximum | 1 | 8 | 8 | 1 | 1 | 1 | 40 |

$q_0$, $\varphi$, $\beta$, $\phi$, $\alpha$, $\rho$ and $na$ are the parameters used in ACS. $\phi$ is explained in the beginning of section III

DPSO = ]0,1]x]-1,8]x]-1,8]x]0,1]x]0,1]x]0,1]x]0,40] $\subset \Re^7$ is the domain of the PSO. We define the fitness value of a given position (point) as the length of the best tour computed by an ACS using the related parameters in the given instance. If comparing two different positions they have the same length value then computing time is considered. We consider better of those parameters that minimize the length of the tour and secondly the time of computing. For computing the fitness of a given position, first integer parameters ($na$, $\varphi$ and $\beta$) are rounded up as shown in Fig. 1, secondly the algorithm runs five trials of the ACS algorithm using the rounded parameters in the TSP-instance and returns the best value obtained from the trials.

Point in the PSO domain (reflects a set of parameters):
(0.1 **2.3 8** 0.5 0.88 0.34 **32.4**)
Modified values to run on ACS:
(0.1 **3 8** 0.5 0.88 0.34 **33** )
Fig. 1. Modification of PSO points. In bold are the modified values

For each particle of the population its initial velocity is set randomly. For half of the population the initial position is set using predefined parameters assuring that for every parameter there will be a particle containing a value covering the full range. The positions of the other half of the initial population is set randomly.

Parameters for the PSO (see Section II) have been set following [13], [14]: $c_1=c_2=2$, $\chi=0.729$, the inertia weight is set initially to 1 and gradually decreasing from 1 to 0.1 (at each PSO iteration $w=0.99w$). Maximum number of iterations is set to 500 due to computing time constraints (for 1000 PSO-iterations more than 1 day of computing time was necessary).

*B. Algorithm*

The algorithm PSO-ACS's pseudo-code is as follows:

Select TSP-Instance.
Initialize particles.
Do 500
   For all the set of particles
     *Position_Fitness* (*PF*)=INFINITE
     Do 5
        Perform a trial ACS with particle parameters.
        If *NewValue < PF*
           *PF=NewValue*
        End if
     End Do
   End for
   Compute $w=0.99w$
   Update Best Parameters Found by each Particle
   Update Best Parameters Found by the Population
   Compute Velocity
   Movement of Particles
End Do
Return the set of parameters related to the best tour length found and the tour length.

The algorithm is based in a PSO framework, where particles are initialized and iteratively are moving though the domain of the set of parameters. The goal of the algorithm is, for a given instance, to compute the tour with lowest length and to compute the set of ACS parameters, among those in DPSO, which gets the best ACS performance. Those final parameters are related with the TSP-instance selected.

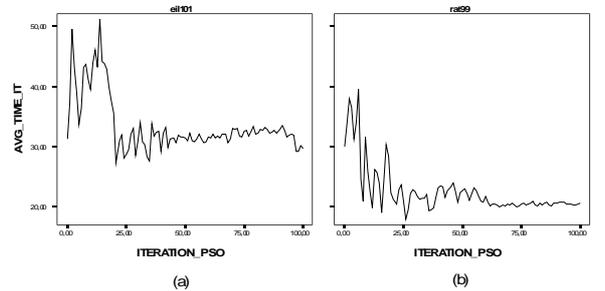

Fig 2. Average time of the swarm at first 100 iterations for instances eil101 (a) and rat99 (b). AVG_TIME is average time of the iteration (given a fixed number of particles) and ITERATION_PSO is the number of the iteration in the PSO-ACS algorithm.

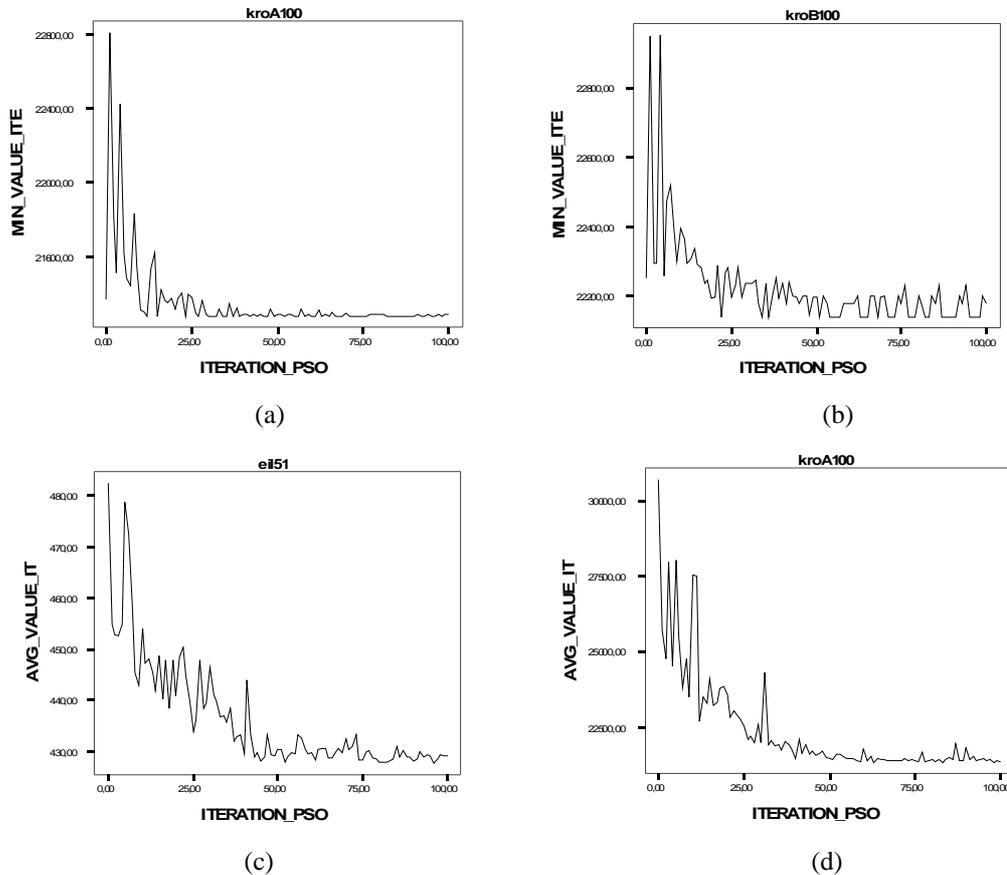

Fig. 3. First 100 iterations of the PSO-ACS algorithm. (a) and (b) are related to the minimum tour obtained at each iteration, (c) and (d) are related to the average of the particles´ fitness values. (a) and (d) are related to the instance kroA100; (b) is related to kroB100; and (c) to eil51. Those are examples of typical behaviors in the 100 first iterations of the algorithm.

TABLE II.
SETS OF PARAMETERS

|  | $\alpha$ | $\beta$ | $\rho$ | $\varphi$ | na | $q_0$ | $\phi$ | Fitness |
|---|---|---|---|---|---|---|---|---|
| P_eil51 | 0.36 | 7 | 0.40 | 1 | 1 | 0.54 | 0.18 | **426** |
| P_eil76 | 0.21 | 5 | 0.40 | 1 | 5 | 0.58 | 0.20 | **538** |
| P_eil101 | 0.71 | 7 | 0.23 | 2 | 7 | 0.78 | 0.12 | **629** |
| P_kroA100 | 0.64 | 4 | 0.24 | 1 | 4 | 0.64 | 0.12 | **21282** |
| P_kroB100 | 0.71 | 1 | 0.08 | 1 | 39 | 0.86 | 0.12 | **22141** |
| P_rat99 | 0.15 | 3 | 0.28 | 1 | 9 | 0.95 | 0.00 | **1211** |
| ACS | 0.10 | 2 | 0.10 | 1 | 10 | 0.9 | a | b |
| ACS_GA | 0.20 | 6 | 0.20 | 1 | 10 | 0.7 | a | b |

"a" have been tested for $\phi = 0.1\ 0.2\ \ldots\ 0.9\ 1$ "b" there is no fitness value related. Values in **bold** mean optimums.

## IV. COMPUTATIONAL RESULTS

Algorithm was coded in C++ in a Pentium M 1.7 with 1GB RAM. Algorithm has been run on six of the most widely used TSP-instances. Computational results are given in four parts: PSO-ACS behavior, PSO-ACS optimum values obtained, best set of parameters and comparison among sets of parameters´ performance.

Computationally each PSO-ACS iteration shows a clear convergence: when the "optimum" (defined by the algorithm) number of ants and $\phi$ are nearly fixed, the computational time is also fixed (see fig. 2). In less than 100 iterations algorithm computes an "optimum" for integer parameters and in 200 iterations there are small differences among the optimum found and the particles´ position for real parameters. In Fig. 3 we can see the evolution of the algorithm in the first 100 iterations. For the average of the fitness of the swarm (in a given iteration), there is a decreasing global tendency, and after iteration 75 we can see the average of the fitness is kept on a fixed range, the size of this range is variable as shown in (c) and (d). For the minimum value obtained by the swarm in a given iteration, computational results show that at the beginning there are increasing and decreasing phases, when the particles are exploring their local optimums and moving also to the global one, but near iteration 100 the minimum is

TABLE III.
AVERAGE PERFORMANCE. BY PARAMETER_SET / INSTANCE

|  | eil51 | eil76 | eil101 | kroA100 | kroB100 | rat99 |
|---|---|---|---|---|---|---|
| P_eil51 | 428.6[ab] | 548.04 | 653.68 | 22147.6 | 22718.92 | 1241.64 |
| P_eil76 | 428.52 | 544.36[ab] | 647.28 | 21600.36 | 22357.24 | 1223.8[b] |
| P_eil101 | 433.92 | 552.88 | 651.2[a] | 24382.36 | 24648.12 | 1251.32 |
| P_kroA100 | 430.28 | 548.52 | 645.6 | 21566.12[a] | 22437.28 | 1234 |
| P_kroB100 | 428.6[b] | 545.68 | 643.52[b] | **21495.4**[b] | 22340.56[ab] | 1228.8 |
| P_rat99 | 429.72 | 553.12 | 662.16 | 21915.24 | 22908.16 | 1242.16[a] |
| PACS-1 | 430.28 | 547.92 | 643.8 | 21625.08 | 22348.16 | 1226.4 |
| PACS-2 | 428.56 | 546.28 | 645.2 | 21640.12 | 22377.28 | 1234.16 |
| PACS-3 | 429.52 | 546.32 | 649.2 | 21641.12 | 22436.6 | 1231.88 |
| PACS-4 | 429.64 | 547.84 | 647.4 | 21621.08 | 22446.4 | 1242.72 |
| PACS-5 | 430.6 | 546.28 | 648.24 | 21589.24 | 22412 | 1237.12 |
| PACS-6 | 430.44 | 546.8 | 647.88 | 21624.2 | 22440.8 | 1241.12 |
| PACS-7 | 430.2 | 549.84 | 649.04 | 21584.28 | 22457.44 | 1238.96 |
| PACS-8 | 430.2 | 549 | 649.76 | 21594.56 | 22578.88 | 1245 |
| PACS-9 | 430.04 | 546.56 | 651.72 | 21590.68 | 22429.6 | 1249.64 |
| ACS_GA-1 | 428.56 | 547.32 | 644.08 | 21674.48 | 22345.48 | 1227.48 |
| ACS_GA-2 | **427.84** | 544.28 | 643.92 | 21669.24 | 22366.84 | 1219.88 |
| ACS_GA-3 | **427.84** | **543.56** | 643.84 | 21594.36 | 22376.16 | 1228.28 |
| ACS_GA-4 | **427.84** | 543.96 | 643.92 | 21514.24 | 22358.48 | **1219.2** |
| ACS_GA-5 | 428.32 | **543.56** | 643.58 | 21589.88 | 22385 | 1222.84 |
| ACS_GA-6 | 428.4 | 544.52 | **642.68** | 21524.84 | 22353.16 | 1224.52 |
| ACS_GA-7 | 428.52 | 548.24 | 644.76 | 21586.16 | 22373.72 | 1225.96 |
| ACS_GA-8 | 428.8 | 544.44 | 643.56 | 21514.8 | 22363.6 | 1222.52 |
| ACS_GA-9 | 428.24 | 542.84 | 643.57 | 21582.36 | 22354.68 | 1223.84 |

"a" expected best value. "b" best value obtained among the new sets of parameters. Values in **bold** are the best values obtained.

TABLE IV.
MINIMUM TOUR LENGTH FOUND. BY PARAMETER_SET / INSTANCE

|  | eil51 | eil76 | eil101 | kroA100 | kroB100 | rat99 |
|---|---|---|---|---|---|---|
| P_eil51 | **426**[a] | 539 | 639 | 21470 | 22404 | 1213 |
| P_eil76 | **426** | **538**[a] | 637 | 21318 | 22237 | **1211** |
| P_eil101 | 428 | 544 | 638[a] | 23764 | 23774 | 1214 |
| P_kroA100 | 427 | **538** | 637 | 21343[a] | 22237 | 1212 |
| P_kroB100 | **426** | 539 | **629** | **21282** | 22199[a] | **1211** |
| P_rat99 | 427 | 547 | 644 | 21560 | 22422 | **1211**[a] |
| PACS-1 | 427 | 540 | 634 | **21282** | 22236 | 1212 |
| PACS-2 | **426** | **538** | 630 | **21282** | 22200 | 1212 |
| PACS-3 | 427 | **538** | 634 | **21282** | 22179 | 1212 |
| PACS-4 | 427 | 540 | 630 | **21282** | 22251 | 1213 |
| PACS-5 | **426** | 540 | 630 | **21282** | 22179 | 1214 |
| PACS-6 | 427 | **538** | 638 | 21305 | 22244 | 1213 |
| PACS-7 | 427 | 540 | 636 | **21282** | 22284 | 1213 |
| PACS-8 | **426** | **538** | 632 | 21292 | 22295 | 1221 |
| PACS-9 | 427 | 540 | 639 | **21282** | 22284 | 1213 |
| ACS_GA-1 | 427 | 539 | 632 | 21330 | 22157[o] | **1211** |
| ACS_GA-2 | **426** | **538** | 630 | 21370 | 22237 | **1211** |
| ACS_GA-3 | **426** | **538** | 631 | 21330 | 22295 | **1211** |
| ACS_GA-4 | **426** | **538** | 634 | 21318 | 22231 | **1211** |
| ACS_GA-5 | **426** | **538** | 630 | 21380 | 22258 | **1211** |
| ACS_GA-6 | **426** | **538** | 630 | 21318 | 22268 | **1211** |
| ACS_GA-7 | **426** | **538** | 630 | 21318 | 22274 | **1211** |
| ACS_GA-8 | **426** | **538** | 630 | 21318 | 22241 | **1211** |
| ACS_GA-9 | **426** | **538** | 630 | 21318 | 22237 | **1211** |

"a" expected best value; "o" best value obtained. Values in **bold** mean optimums. For the PACS and ACS_GA set of parameters PACS-i means $\phi=0.i$ maintained as in (a) or frequently visited as in (b). This fast convergence can be an advantage as well as a drawback because it can lead to a fast non-desirable convergence. We set the reasons of this fast convergence in the PSO framework used and mainly in the method for evaluating a set of parameters: in a stochastic algorithm there is the probability that a bad set of parameters could perform well, if all the particles move into this area and the number of iterations in this area increases leading to probably good solutions that will cause the algorithm to remain in this non-optimal area.

Table II shows the optimum set of parameters found running PSO-ACS on each one of the instances selected; P_eil51 is the set of parameters obtained running PSO-ACS on the instance eil51, similarly for P_eil76, P_eil101..and so on.. PACS and ACS_GA are the set of parameters proposed by [6] and [16] respectively, considering several values of $\phi$ (parameter that defines the neighbourhood). "Fitness" is the minimal length tour obtained by the PSO-ACS algorithm in the related instance. There is no clear cut rule to define the optimal parameters even if the recommended ranges can be guessed, also $\beta$ is normally bigger than $\varphi$; $q_0$ is always bigger than 0.5; $\phi$, percentage of vertices in the neighbour, is normally less or equal to the 20%; the number of ants are normally between 1 and 10. But we cannot set a rule because the set of parameters for P_kroB100 is not following those guidelines and is performing better than the rest. In all the instances the optimum is found but this is not relevant because of the enormous quantity of time expended (sometimes more than half a day) and because of the number of times an ACS algorithm is run on a given instance.

We have run each one of the obtained set of parameters in all of the previous selected instances, to check the efficiency of the parameters. As such ACS is performed with 25 trials of 2500 iterations each.

Table III and IV show the performance, based on the average tour length (trial-based) and the minimum tour length found respectively, for each of the sets of parameters on each of the instances. We expect that the parameters related to an instance will compute the best results on this instance. Comparing this hypothesis on Tables III and IV, computational results show that this is in fact false. Instead, working only in this reduced set of instances, the set of parameters for P_kroB100 can be considered as the best overall one which is performing optimally in most of the

TABLE V.
TIME COMPARISON BEST SETS OF PARAMETERS

|  | Avg. FI | Avg. Time |
|---|---|---|
| P_kroB100 | 4,53 | 241,29 |
| Best of PACS | **1,35** | **69,53** |
| Best of GA-ACS | 3,16 | 162,89 |

"Avg. FI" is the average time (considering all refereed instances) used for finding the best solution proposed by the algorithm, considering the time from the trial this solution was computed. "Avg. Time" is the average computational time of the algorithm running in all the proposed instances. In bold are the best times. For PACS and GA_ACS lowest times for all possible $\Phi$ values are used: are named respectively "Best of PACS" and "Best of GA-ACS".

instances and with an efficient average value.

Comparing the sets of parameters of PACS and ACS_GA, one can observe that PACS is performing better in bigger instances and ACS_GA in smaller ones. Further, ACS_GA performs better on the average.

Computational times of PACS, ACS_GA and P_kroB100 set of parameters are compared in Table V. PACS perform best in both measures and P_kroA100 is the worst. The set of parameters of P_kroB100 is performing better than ACS and ACS_GA in most of the instances, but it has a greater computational time.

## V. CONCLUSIONS

Computational results seem to show that there is no uniquely optimal set of ACS parameters yielding best quality solutions in all the TSP instances. Nevertheless the PSO-ACS has been able to find a set of ACS parameters that work optimally for a majority of instances unlike others known so far.

PSO-ACS algorithm works well across different instances because it adapts itself to the instance characteristics. But it has a high computational overhead. A future work will try to modify the algorithm framework to reduce this cost.

PSO-ACS also has a fast convergence that can lead to a bad set of parameters. This may be due to two reasons: first is the specific PSO framework used, and in modifying it we expect to obtain better results. Secondly the way the sets of parameters are evaluated may have to be reviewed as a bad set of parameters could lead to a non-desired convergence.

## ACKNOWLEDGMENT

D. Gómez-Cabrero thanks the University of Colombo School of Computing for the support given while he was in Sri Lanka. The contribution by D. Gómez-Cabrero has been partially supported by the AVCiT of the Generalitat Valenciana (Ref: GRUPOS03/174).